Kernel Groups and nontrivial Galois module structure of imaginary quadratic fields

Daniel R. Replogle

Department Mathematics and Computer Science, College of Saint Elizabeth

Morristown, NJ 07960

Abtract: Let $K$ be an algebraic number field with ring of integers $\mathcal{O}_K$, $p > 2$ be a rational prime and $G$ be the cyclic group of order $p$. Let $\Lambda$ denote the order $\mathcal{O}_K[G]$. Let $Cl(\Lambda)$ denote the locally free class group of $\Lambda$ and $D(\Lambda)$ the kernel group, the subgroup of $Cl(\Lambda)$ consisting of classes that become trivial upon extension of scalars to the maximal order. If $p$ is unramified in $K$, then $D(\Lambda) = T(\Lambda)$, where $T(\Lambda)$ is the Swan subgroup of $Cl(\Lambda)$. This yields upper and lower bounds for $D(\Lambda)$. Let $R(\Lambda)$ denote the subgroup of $Cl(\Lambda)$ consisting of those classes realizable as rings of integers, $\mathcal{O}_L$, where $L/K$ is a tame Galois extension with Galois group $Gal(L/K) \cong G$. We show under the hypotheses above that $T(\Lambda)^{(p-1)/2} \subseteq R(\Lambda) \cap D(\Lambda) \subseteq T(\Lambda)$, which yields conditions for when $T(\Lambda) = R(\Lambda) \cap D(\Lambda)$ and bounds on $R(\Lambda) \cap D(\Lambda)$. We carry out the computation for $K = \mathbb{Q}(\sqrt{-d})$, $d > 0$, $d \neq 1$ or $3$. In this way we exhibit primes $p$ for which these fields have tame Galois field extensions of degree $p$ with nontrivial Galois module structure.

**Correspondence Address:**
  Daniel R Replogle
  Department of Mathematics and Computer Science
  College of Saint Elizabeth
  2 Convent Road
  Morristown, NJ 07960

**Email: dreplogle@liza.st-elizabeth.edu**

Section 1: Introduction and Subgroups of $Cl(\Lambda)$

Let $K$ be an algebraic number field and denote its ring of algebraic integers by $\mathcal{O}_K$. Let $G$ be a finite abelian group of order $n$. Let $\Lambda$ denote the order $\mathcal{O}_K[G]$ in the group algebra $K[G]$. The class group of stable isomorphism classes of locally free $\Lambda$-modules is denoted by $Cl(\Lambda)$. The kernel group, $D(\Lambda)$, is the subgroup of $Cl(\Lambda)$ consisting of those classes that become trivial upon extension of scalars to the maximal order. Let $\Sigma = \sum_{g \in G} g$ then for each $r \in \mathcal{O}_K$ so that $r$ and $n$ are relatively prime define the Swan module $\langle r, \Sigma \rangle$ by $\langle r, \Sigma \rangle = r\Lambda + \Lambda\Sigma$. Swan modules are locally free rank one $\Lambda$-ideals and hence determine classes in $Cl(\Lambda)$ [15]. The set of classes of Swan modules is the Swan subgroup, denoted $T(\Lambda)$, a subgroup of $D(\Lambda)$ first extensively studied by Ullom [15]. (The fact that $T(\Lambda)$ is a subgroup of $D(\Lambda)$ follows from the Mayer-Vietoris sequence of Reiner-Ullom discussed in section 2). Throughout this article let $p > 2$ be prime, $\zeta_p$ be a primitive $p$-th root of unity, and $\mathbb{Q}$ the field of rational numbers. Let $C_s$ denote the cyclic group of order $s$. In our results we will require that $p > 2$ is a rational prime so that $p$ is unramified in $\mathcal{O}_K$. We will express this as saying $p > 2$ is unramified in $K$. Recall two number fields $E$ and $F$ are arithmetically disjoint if no rational primes ramifiy in both $\mathcal{O}_E$ and $\mathcal{O}_F$ and in this case for $K = EF$, the compositum, we have $\mathcal{O}_K = \mathcal{O}_E\mathcal{O}_F$ (see [1] for example).

The first main result of this article is the following result of section 2.

**Theorem 2.1.** *Let $p > 2$ be unramified in a number field $K$ and let $G \cong C_p$ then $T(\Lambda) = D(\Lambda)$.*

Since $T(\Lambda)$ is easily seen to be trivial for cyclic groups $G$ when $K = \mathbb{Q}$, Theorem 2.1 implies that $D(\Lambda)$ is trivial when $K = \mathbb{Q}$, a well-known theorem of Rim [13].

Let $R(\Lambda)$ denote the subgroup of $Cl(\Lambda)$ consisting of those classes realizable as rings of integers $\mathcal{O}_L$ where $L/K$ is a tame abelian Galois extension of number fields with abelian Galois group $Gal(L/K) \cong G$ [9]. $R(\Lambda)$ is described explicitly as a subgroup of $Cl(\Lambda)$ for all $p$-elementary abelian groups in [8]. Note that showing $R(\Lambda)$ is nontrivial proves the existence of a tame Galois field extension $L/K$ with $Gal(L/K) \cong G$ so that $\mathcal{O}_L$ is not a free $\mathcal{O}_K[G]$-module. In the language of [6] this shows $K$ has nontrivial Galois module structure for $G$. (Equivalently, we may say in this case there exists an $L$ so that $L/K$ does not have a normal integral basis. We note from [8] it follows there are infinitely many such fields for each nontrivial class in $R(\Lambda)$). Using McCulloh's description of $R(\Lambda)$ from [8] and the relationship between $R(\Lambda)$, $T(\Lambda)$, and $D(\Lambda)$ from [6] we obtain the second main result of this article (which will be stated in a stronger form in Section 4).

**Theorem 4.2.** *Let $G \cong C_p$ and $p > 2$ be unramified in $K$. If the exponent of $T(\Lambda)$ is relatively prime with $(p-1)/2$ then $T(\Lambda) = R(\Lambda) \cap D(\Lambda)$.*

We note the proof of [6, Theorem 2] shows for any algebraic number field $K \neq \mathbb{Q}$ there are in fact infinitely many primes $p$ so that for each there is a tame Galois field extension $L$ of $K$ of degree $p$ without a normal integral basis (i.e. for which $L/K$ has a nontrivial Galois module structure).

We now outline the structure of this article. In Section 2 we introduce the exact Mayer-Vietoris sequence of Reiner-Ullom which yields a convenient description of the Swan subgroup $T(\Lambda)$. We use this to prove Theorem 2.1 and recover Rim's Theorem from Theorem 2.1. In Section 3 we combine Theorem 2.1 with upper and lower bounds for Swan subgroups to obtain upper and lower bounds for kernel groups. We carry out the computation for imaginary quadratic fields $K = \mathbb{Q}(\sqrt{-d})$ when $d > 0$ and $d \neq 1$ or $3$. In Section 4 we explore the implications this work has to Galois module structure problems. Specifically we prove Theorem 4.2 and derive from it our concluding result our last main theorem:

**Theorem 4.3.** *Let $K = \mathbb{Q}(\sqrt{-d})$ where $d > 0$ and $d \neq 1$ or 3. Then*

(a) *For $p$ inert in $K/\mathbb{Q}$ we have $C_{\frac{p+1}{2}} \subseteq R(\Lambda) \cap D(\Lambda) \subseteq C_{p+1}$.*

(b) *For $p$ split in $K/\mathbb{Q}$ we have $R(\Lambda) \cap D(\Lambda)$ is a homomorphic image of $C_{p-1}$.*

(c) *For $p$ ramified in $K/\mathbb{Q}$ we have $T(\Lambda) \cong C_p$, hence $C_p \subseteq R(\Lambda) \cap D(\Lambda)$.*

SECTION 2: SWAN SUBGROUPS AND A GENERALIZATION OF RIM'S THEOREM

The key properties of the Swan subgroup used in this article are the exact Mayer-Vietoris sequence of Reiner-Ullom and the convenient description of the Swan subgroup it yields. (Here and throughout the rest of the text we will denote $\mathcal{O}_K$ by $\mathcal{O}$ when no confussion can result). Let $\Gamma = \Lambda/(\Sigma), \overline{\mathcal{O}} = \mathcal{O}/p\mathcal{O}$, $\phi$ and $\overline{\phi}$ denote the canonical quotient maps, $\epsilon$ denote the augmentation map, and $\overline{\epsilon}$ denote the map induced by the augmentation map. Consider the fiber product:

$$\begin{array}{ccc} \Lambda & \xrightarrow{\phi} & \Gamma \\ \epsilon \downarrow & & \downarrow \overline{\epsilon} \\ \mathcal{O} & \xrightarrow{\overline{\phi}} & \overline{\mathcal{O}}. \end{array}$$

The result in [10] applied to the case when $G$ is abelian (or more generally, when the group algebra $K[G]$ satisfies the "Eichler condition"– see [10] or [3]) is: There is an exact Mayer-Vietoris sequence:

$$1 \longrightarrow \Lambda^* \longrightarrow \mathcal{O}^* \times \Gamma^* \xrightarrow{h} \overline{\mathcal{O}}^* \xrightarrow{\delta} D(\Lambda) \xrightarrow{\phi} D(\Gamma) \oplus D(\mathcal{O}) \longrightarrow 0,$$

where for any ring $S$ we denote its group of units by $S^*$. From [15] we have that the image of $s \in \overline{\mathcal{O}}^*$ under $\delta$ is $[s, \Sigma]$, the class of the Swan module $\langle s, \Sigma \rangle$, and hence $T(\Lambda) = Im(\delta)$; therefore, $T(\Lambda)$ is a subgroup of $D(\Lambda)$ and $T(\Lambda) \cong \overline{\mathcal{O}}^*/h(\mathcal{O}^* \times \Gamma^*)$. We note that the map $h$ is given by $(u, v) \mapsto \overline{u} \cdot \overline{v}^{-1} = \phi(u)\overline{\epsilon}(v)^{-1}$.

We now prove our first main result which gives condtions when $D(\Lambda)$ and $T(\Lambda)$ coincide.

**Theorem 2.1.** *Let $p > 2$ be unramified in a number field $K$ and let $G \cong C_p$ then $T(\Lambda) = D(\Lambda)$.*

*Proof.* Since $p$ does not ramify in $\mathcal{O}$, $\mathcal{O}/p\mathcal{O} \cong \mathcal{O}[\zeta_p]/(1 - \zeta_p)$. We also have that $\mathcal{O}[\zeta_p] \cong \Lambda/(\Sigma)$. Let $x$ be a fixed generator of $G = C_p$. Let $f : \Lambda = \mathcal{O}[G] \longrightarrow \mathcal{O}[\zeta_p]$ be induced by $x \mapsto \zeta_p$ and $j : \Lambda = \mathcal{O}[G] \longrightarrow \mathcal{O}$ be induced by $x \mapsto 1$. Then we have the fiber product:

$$\begin{array}{ccc} \Lambda & \xrightarrow{f} & \mathcal{O}[\zeta_p] \cong \Lambda/\mathcal{O}\Sigma = \Gamma \\ j \downarrow & & \downarrow \\ \mathcal{O} & \longrightarrow & \mathcal{O}[\zeta_p]/(1 - \zeta_p) \cong \overline{\mathcal{O}} \end{array}$$

This gives rise to the exact Mayer-Vietoris sequence:

$$1 \longrightarrow \Lambda^* \longrightarrow \mathcal{O}^* \times \Gamma^* \xrightarrow{h} \overline{\mathcal{O}}^* \xrightarrow{\delta} D(\Lambda) \longrightarrow D(\mathcal{O}) \oplus D(\Gamma) \longrightarrow 0.$$

Since $K$ and $\mathbb{Q}(\zeta_p)$ are arithmetically disjoint, $\Gamma \cong \mathcal{O}[\zeta_p]$, a maximal order in $K\mathbb{Q}(\zeta_p)$, whence $D(\Gamma) \cong \{1\}$. Similarly, $\mathcal{O}$ is a maximal order in $K$ so $D(\mathcal{O}) \cong \{1\}$. So our Mayer-Vietoris sequence becomes:

$$1 \longrightarrow \Lambda^* \longrightarrow \mathcal{O}^* \times \Gamma^* \xrightarrow{h} \overline{\mathcal{O}}^* \xrightarrow{\delta} D(\Lambda) \longrightarrow 0.$$

Hence $\delta$ is surjective implying $D(\Lambda) = T(\Lambda)$. □

**Corollary 2.2 (Rim's Theorem [13]).** *For each odd prime $p$ and $G \cong C_p$, $D(\mathbb{Z}[G]) \cong \{1\}$.*

*Proof.* By Theorem 2.1 it suffices to show $h$ is surjective. The argument is standard and is essentially the proof that $T(\mathbb{Z}[G])$ is trivial whenever $G$ is cyclic, (see [Section 53, 3] for example). To simplify notation let $\pi = 1 - \zeta_p$ and $S = \mathbb{Z}[\zeta_p]/(\pi)$. Each element of $S^*$ is of the form $\overline{\eta} = n + \pi\mathbb{Z}[\zeta_p]$ for some $n \in \mathbb{Z}$ relatively prime to $p$. Since the cyclotomic unit $u = (1 - \zeta_p{}^n)/(1 - \zeta_p) = [1 - (1 - \pi)^n]/\pi$ of $\mathbb{Z}[\zeta_p]$ is congruent to $n \mod \pi$, $h[(u, 1)] = \overline{\eta}$. Hence, $h$ is surjective. □

SECTION 3: UPPER AND LOWER BOUNDS FOR KERNEL GROUPS

In this section we outline proofs of upper and lower bounds for the Swan subgroup, adapted from [6], [11], and [12]. Combining these bounds with Theorem 2.1 yields upper and lower bounds on the kernel group, $D(\Lambda)$ when $p > 2$ is unramified in $K$ and $G \cong C_p$.

We first introduce some notation and state our upper and lower bounds as Proposition 3.1. We then give the arguments.

Let $\overline{\mathcal{O}} = \mathcal{O}/p\mathcal{O}$. As above $\mathcal{O}^*$ denotes the units of $\mathcal{O}$. We denote by $Im(\mathcal{O}^*)$ the image of $\mathcal{O}^*$ via the map $s \mapsto \overline{s}$ where $s \in \mathcal{O}^*$ and $\overline{s}$ is the image of $s$ under the canonical homomorphism $can : \mathcal{O} \to \mathcal{O}/p\mathcal{O}$. Let $\Delta = Aut(G)$, the group of automorphisms of $G$. Let $V_p = \overline{\mathcal{O}}^*/Im(\mathcal{O}^*)$, and for an abelian group $H$ we denote by $(H)_k$ its $k$-torsion subgroup.

**Proposition 3.1.** *If $p > 2$ is unramified in $K$ and $G \cong C_p$, then*

$$V_p^{p-1} \leq D(\Lambda) \leq \overline{\mathcal{O}}^*/(Im(\mathcal{O}^*)(\mathbb{Z}/p\mathbb{Z})^*),$$

*where "$\leq$" is taken to mean "isomorphic to a subgroup of."*

*Proof.* In Lemma 3.2 below we will show the existence of a surjective map $T \longrightarrow V_p^{p-1}$. Hence we have $V_p^{p-1} \leq T(\Lambda)$. In Lemma 3.4 we obtain that $T(\Lambda) \leq \overline{\mathcal{O}}^*/Im(\mathcal{O}^*)(\mathbb{Z}/p\mathbb{Z})^*$. Hence we have

$$V_p^{p-1} \leq T(\Lambda) \leq \overline{\mathcal{O}}^*/Im((\mathcal{O}^*)(\mathbb{Z}/p\mathbb{Z})^*).$$

The proposition now follows from Theorem 2.1. □

**Lemma 3.2 ([Theorem 12, 12]).** *There is a surjective map $T(\Lambda) \longrightarrow V_p/(V_p)_{p-1} \cong V_p^{p-1}$.*

*Sketch of Proof.* We have $T(\Lambda) \cong (\overline{\mathcal{O}}^*/Im(\mathcal{O}^*))/(Im(\mathcal{O}^* \times \Gamma^*)/Im(\mathcal{O}^*))$ under the map $h : \mathcal{O}^* \times \Gamma^* \longrightarrow \overline{\mathcal{O}}^*$ given by $(s, \gamma) \mapsto \overline{s}\,\overline{\gamma}^{-1}$. It suffices to show that

$$h(\mathcal{O}^* \times \Gamma^*)/Im(\mathcal{O}^*) \subseteq \{x \in \overline{\mathcal{O}}^*/Im(\mathcal{O}^*) : x^{p-1} = 1\}$$

To see this it suffices to show that for any $\gamma \in \Gamma^*$ that $\bar{\epsilon}(\gamma)^{p-1} \in Im(\mathcal{O}^*)$.

One shows there is an isomorphism $\mathcal{O}^* \cong (\Gamma^\Delta)^*$. Then let $N$ be the norm map $N : \Gamma^* \longrightarrow (\Gamma^\Delta)^* \cong \mathcal{O}^*$, defined by $N(\gamma) = \prod_{\delta \in \Delta} \gamma^\delta$. One then shows that the diagram

$$\begin{array}{ccc} \Gamma^* & \xrightarrow{N} & (\Gamma^\Delta)^* \cong \mathcal{O}^* \\ \downarrow \bar{\epsilon} & \downarrow \bar{\epsilon} & \downarrow can \\ \overline{\mathcal{O}}^* & \xrightarrow{(\ )^{p-1}} & \overline{\mathcal{O}}^* = \overline{\mathcal{O}}^* \end{array}$$

commutes, completing the proof. □

**Remarks 3.3.** *1. For details on the above see [12, Lemmata 13 and 14] or [6, Theorem 5]. 2. The proof in [6] gives a generalization to the case $G$ is $p$-elementary abelian. 3. The proofs of Lemmata 3.2 and 3.4 implicitly use the finiteness of the classgroup to get the "isomorphic to a subgroup of" claim.*

**Lemma 3.4.** *If $G \cong C_n$ and $K$ is a number field then $T(\Lambda) \leq \overline{\mathcal{O}}^*/((\mathbb{Z}/n\mathbb{Z})^* Im(\mathcal{O}^*))$*

*Proof.* Let $G = C_n$ where $n \geq 3$ (for $n = 1$ or $2$ there is nothing to prove). It is well known that when $G$ is cyclic $T(\mathbb{Z}[G])$ is trivial. Therefore by extending scalars from $\mathbb{Z}$ to $\mathcal{O}$ we can conclude the Swan class $[s, \Sigma] = 0 \in T(\Lambda)$ whenever $s \in \mathbb{Z}$. In particular $[r\Lambda + \Sigma\Lambda] = 0$ for all $r \in \mathbb{Z}$ with $(r, n) = 1$. Hence the kernel of the surjective map $\rho : \overline{\mathcal{O}}^*/Im(\mathcal{O}^*) \to T(\Lambda)$ contains $(\mathbb{Z}/n\mathbb{Z})^*$.

Therefore $\rho$ yields a map from

$$(\overline{\mathcal{O}}^*/Im(\mathcal{O})^*)/(((\mathbb{Z}/n\mathbb{Z})^* Im(\mathcal{O}^*))/Im(\mathcal{O}^*)) \cong \overline{\mathcal{O}}^*/((\mathbb{Z}/n\mathbb{Z})^* Im(\mathcal{O}^*))$$

onto $T(\Lambda)$. □

We now consider imaginary quadratic extensions of $\mathbb{Q}$. Let $K = \mathbb{Q}(\sqrt{-d})$, where $d$ is positive and for convenience we assume $d \neq 1$ or $3$ (but see Section 4 and [12]). We wish to obtain information about $D(\Lambda)$ when $G \cong C_p$ and $p > 2$ is unramified in $K$. The result we prove is given next.

**Proposition 3.5.** *Let $K = \mathbb{Q}(\sqrt{-d})$ where $d > 0$ is square free and coprime to $p$, and $d \neq 1$ or $3$. Let $G$ be cyclic of order $p$. Then $D(\Lambda)$ is a homomorphic image of $C_{p-1}$ if $p$ splits completely in $K$, and $D(\Lambda)$ is isomorphic to $C_{p+1}$ or $C_{(p+1)/2}$ if $p$ is inert in $K$.*

We first compute $\overline{\mathcal{O}}/Im(\mathcal{O}^*)$ for imaginary quadratic fields. Then we compute the upper and lower bounds of $T(\Lambda)$ for $d$ as above. The result then immediately follows from Proposition 3.1. (For the case $G$ is of order 2 analogous results hold, see [14]. The finite field of $p^n$ elements will be denoted $\mathbf{F}_{p^n}$.

*Proof.* If $p$ is inert, then $\mathcal{O}/p\mathcal{O}$ is a field and has degree 2 over $\mathbb{Z}/p\mathbb{Z}$. Thus $\overline{\mathcal{O}} = \mathcal{O}/p\mathcal{O} \cong \mathbf{F}_{p^2}$ whence $\overline{\mathcal{O}}^* \cong C_{p^2-1}$. If $p$ splits then $(p)$ factors into two distinct prime ideals in $\mathcal{O}$, say $(p_1)$ and $(p_2)$. Then by the Chinese remainder theorem we get $\mathcal{O}/p\mathcal{O} \cong \mathcal{O}/p_1\mathcal{O} \times \mathcal{O}/p_2\mathcal{O}$, whence $\mathcal{O}/p\mathcal{O} \cong \mathbf{F}_p \times \mathbf{F}_p$ and $\overline{\mathcal{O}}^* \cong C_{p-1} \times C_{p-1}$.

Next note that $Im(\mathcal{O}^*) \subset (\mathbb{Z}/p\mathbb{Z})^*$, so the upper bound of Proposition 3.2 is $\overline{\mathcal{O}}^*/((\mathbb{Z}/p\mathbb{Z})^*)$. If $p$ splits completely we have $\overline{\mathcal{O}}^* \cong C_{p-1} \times C_{p-1}$, which forces the lower bound $V_p^{p-1}$ to be trivial. For the upper bound we have $\overline{\mathcal{O}}^*/(\mathbb{Z}/p\mathbb{Z})^* \cong C_{p-1}$. Therefore $T(\Lambda)$ is a homomorphic image of $C_{p-1}$.

If $p$ is inert we have $\overline{\mathcal{O}}^* \cong C_{p^2-1}$. Therefore our lower bound is given by $C_{(p+1)/2}$. For the upper bound we have $\overline{\mathcal{O}}^*/(\mathbb{Z}/p\mathbb{Z})^* \cong C_{p+1}$. Therefore we obtain that $T(\Lambda)$ is isomorphic with one of these two groups, i.e. $T(\Lambda) \cong C_{p+1}$ or $C_{(p+1)/2}$. □

We note two facts of interest. First one can compute $T(\Lambda)$ when $p$ ramifies exactly as above, however then Theorem 2.1 is not applicable. We include this for completeness as our last proposition of this section. Second, by comparison, Fröhlich, Reiner and Ullom [4], [5], [10] proved that $D(\mathbb{Z}[G])$ is trivial if $G \cong C_p$, and for $G$ an abelian $p$-group of order $p^n$ with $n > 1$ that $D(\mathbb{Z}[G])$ is a $p$-group. Hence we see even in the simplest case, that of an imaginary quadratic field and $G \cong C_p$, their results do not generalize to rings of algebraic integers.

**Proposition 3.6.** *Let $K = \mathbb{Q}(\sqrt{-d})$ where $d > 0$ and $d \neq 1$ or $3$. Assume the prime $p > 2$ divides $d$ and let $G \cong C_p$. Then $T(\Lambda) \cong C_p$.*

*Proof.* Since $p$ divides $d$ we have $p$ ramifies. Thus $p\mathcal{O}_K = \mathcal{P}^2$ for some prime ideal $\mathcal{P}$ in $\mathcal{O}_K$. Moreover we can write as an equation of ideals $p\mathcal{O}_K = (p, \mu)$ where $\mathcal{O}_K$ has a maximal ideal of the form $p\mathcal{O}_K + \mu\mathcal{O}_K$ and $\mathcal{O}_K/p\mathcal{O}_K \cong \mathbf{F}_p[\mu]$ where $\mu^2 = 0$. Hence we may write $\mathcal{O}_K/p\mathcal{O}_K = \{a + b\mu \ : \ a, b \in \mathbf{F}_p\}$. The units in this ring are seen to be of the form $U = \{a + b\mu \ : \ a \neq 0\}$. The sets $\{1 + b\mu : b \in \mathbf{F}_p\}$ and $\{a + 0\mu : a \neq 0\}$ form subgroups whose union is all of $U$ and whose intersection is the identity and are isomoprhic to $C_p$ and $C_{p-1}$, respectively.

Thus we have $(\mathcal{O}_K/p\mathcal{O}_K)^* \cong C_{p-1} \times C_p$. Now as $p > 2$ is odd one has $Im(\mathcal{O}_K^*) \cong C_2$. Again as $p$ is odd $Im(\mathcal{O}_K^*)$ must embed in the $C_{p-1}$ factor. Hence $\overline{\mathcal{O}_K}^*/Im(\mathcal{O}_K^*) \cong C_{\frac{p-1}{2}} \times C_p$. The upper and lower bounds are thus both $C_p$. This gives $T(\Lambda) \cong C_p$. However as $p$ ramifies we need not have $D(\Lambda) = T(\Lambda)$ hence we obtain no upper bound on $D(\Lambda)$. □

SECTION 4: NONTRIVIAL $p$-EXTENSIONS OF IMAGINARY QUADRATIC FIELDS

Let $G$ be the additive group of $\mathbf{F}_{p^n}$. Let $L$ range over all Galois extensions of $K$ with Galois group isomorphic to $G$. It is well known that $\mathcal{O}_L$ is a locally free $\Lambda = \mathcal{O}_K[G]$-module if and only if the extension $L/K$ is tame (i.e. at most tamely ramified). Therefore, when the extension $L/K$ is tame $\mathcal{O}_L$ yields a "Galois module class" in $Cl(\Lambda)$. Following McCulloh ([8]) denote by $R(\Lambda)$ the set of all classes in $Cl(\Lambda)$ realizable as rings of integers of tame Galois extensions of $K$ with Galois group $G$.

Denote by $C$ the group of nonzero elements of $\mathbf{F}_{p^n}$. Since $C \subset Aut(G)$, $Cl(\Lambda)$ is a $\mathbb{Z}[C]$-module. For each $\delta \in C$ denote by $t(\delta)$ the least nonnegative residue (mod $p$) of $Tr(\delta)$ where $Tr: \mathbf{F}_{p^n} \longrightarrow \mathbf{F}_p$ is the trace. Let $\theta = \sum_{\delta \in C} t(\delta)\delta^{-1} \in \mathbb{Z}[C]$, and let $\mathcal{J} = \mathbb{Z}[C](\theta/p) \cap \mathbb{Z}[C]$, the Stickelberger ideal of $\mathbb{Z}[C]$. Then the main theorem in [8] is:

**Theorem (McCulloh, [8]).** *If $G$ is elementary abelian of order $p^n$, then $R(\Lambda) \cong Cl^0(\Lambda)^{\mathcal{J}}$, where $Cl^0(\Lambda)$ is the kernel of the map $Cl(\Lambda) \longrightarrow Cl(\mathcal{O}_K)$ induced by the augmentation $\mathcal{O}_K[G] \to \mathcal{O}_K$.*

Since $T(\Lambda)$ and $D(\Lambda)$ are $\mathbb{Z}[C]$-submodules of $Cl^0(\Lambda)$, McCulloh's theorem implies the following proposition.

**Proposition 4.1 ([Corollary 7, 6]).** *For $G$ elementary abelian of order $p^n \neq 2$, $T^{p^{n-1}(p-1)/2}(\Lambda) = T(\Lambda)^{\mathcal{J}} \leq R(\Lambda) \cap D(\Lambda)$.*

This proposition follows from two facts. First, $C$ acts trivially on a Swan class. Second, for $\varepsilon: \mathbb{Z}[C] \longrightarrow \mathbb{Z}$ the augmentation, if $p^n \neq 2$ then $\varepsilon(\mathcal{J}) = (p^{n-1}(p-$

1)/2)$\mathbb{Z}$. For details see [6].

This last result yields some information about $R(\Lambda) \cap D(\Lambda)$, namely we prove Theorem 4.2 which we now state in a stronger form.

**Theorem 4.2.** *Suppose $G \cong C_p$ and $p > 2$ is unramified in $K$. Then:*

(a) $V_p^{(p-1)^2/2} \leq R(\Lambda) \cap D(\Lambda) \leq \overline{\mathcal{O}}^* / (Im(\mathcal{O}^*)(\mathbb{Z}/p\mathbb{Z})^*).$

(b) *If the exponent of $T(\Lambda)$ is coprime to $(p-1)/2$, then $R(\Lambda) \cap D(\Lambda) = T(\Lambda)$.*

*Proof.* From Theorem 2.1 we have $T(\Lambda) = D(\Lambda)$. From the Proposition 4.1 we have, as $n = 1$, $T^{(p-1)/2}(\Lambda) \subseteq R(\Lambda) \cap D(\Lambda)$. Thus it follows $T^{(p-1)/2}(\Lambda) \subseteq R(\Lambda) \cap D(\Lambda) \subseteq D(\Lambda) = T(\Lambda)$. The lower and upper bounds of (a) now follow from Lemmata 3.2 and 3.4. If the exponent of $T(\Lambda)$ is relatively prime to $(p-1)/2$, then $T^{(p-1)/2}(\Lambda) = T(\Lambda)$, giving (b). $\square$

We note this result generalizes [Corollary 7, 6] by also providing an upper bound and a condition for equality. We now give the application of Theorem 4.2 announced in the introduction.

**Theorem 4.3.** *Let $K = \mathbb{Q}(\sqrt{-d})$, $d$ square-free and $d \neq 1$ or $3$.*

(a) *If $p$ is inert in $\mathcal{O}_K$ then $C_{(p+1)/2} \leq R(\Lambda) \cap D(\Lambda) \leq C_{p+1}$.*

(b) *If $p$ splits in $\mathcal{O}_K$ then $R(\Lambda) \cap D(\Lambda)$ is a homomorphic image of $C_{p-1}$.*

(c) *If $p$ ramifies in $K/\mathbb{Q}$ then $T(\Lambda) \cong C_p$. Hence $C_p \leq R(\Lambda) \cap D(\Lambda)$.*

*Proof.* This follows from Theorem 4.1 (a) and the computations in Section 3. For instance for (a) we have: $V_p \cong C_{(p^2-1)/2}$ and $\overline{\mathcal{O}}^* / Im(\mathcal{O}(\mathbb{Z}/p\mathbb{Z})^*) \cong C_{p+1}$.

For (c) see Proposition 3.6 and note $C_p$ is not of exponent dividing $p - 1$. Also note, Theorem 2.1 does not apply hence we do not obtain an upper bound on $R(\Lambda) \cap D(\Lambda)$ and only obtain a lower bound on $D(\Lambda)$. $\square$

We close by noting the cases $d = 1$ or $d = 3$ were only omitted as there are more units to consider. That is for $K = \mathbb{Q}(\sqrt{-1})$ the units in $\mathcal{O}_K$ are given by $\{\pm 1, \pm i\}$. Determining how these embed in $\mathcal{O}_K/p\mathcal{O}_K$ is not difficult. Likewise, for $K = \mathbb{Q}(\sqrt{-3})$ the units are $\{\pm 1, \pm \omega, \pm \omega^2\}$ where $\omega$ is a primitive sixth root of unity. The problem is of course one has to examine many cases. Moreover, generalizing these results to biquadratic fields (compositums of the types of fields considered) presents only one other difficulty, namely computing the embedding of the fundamental unit. More precisely, let $K_1 = \mathbb{Q}(\sqrt{-d_1})$ and $K_2 = \mathbb{Q}(\sqrt{-d_2})$, where $d_1$, $d_2 > 0$, and so that the discriminants are relatively prime. Then $K_1$ and $K_2$ are arithmetically disjoint, so for $K = K_1 K_2$, the compositums, we have $\mathcal{O}_K = \mathcal{O}_{K_1} \mathcal{O}_{K_2}$. Moreover, the only other intermediate field of $K/\mathbb{Q}$ is $K_3 = \mathbb{Q}(\sqrt{d_1 d_2})$. It is easy to see using Dirichlet's unit theorem that $\mathcal{O}_K^*$ contains the roots of unity in $K$ and a set of units generated by one fundamental unit. Hence one may do the computations required to find bounds on the Swan subgroup. The problem is of course determining the embedding of this fundamental unit in the appropriate quotient and then computing the quotient for each possible case.

It might be of interest then to consider one particular case $K = \mathbb{Q}(\sqrt{-1}, \sqrt{-3})$. Note that $\mathbb{Q}(\sqrt{-1}, \sqrt{-3}) = \mathbb{Q}(\zeta_{12})$ where $\zeta_{12}$ is a primitve 12th root of unity. In [2] several primes are exhibited for which this field has a nontrivial Galois module structure. Hence using the methods contained in here for this field and those primes might be of interest as the Swan subgroups are not explicitly computed in [2]. This project could lead to other interesting ideas and might be purused in a furture note.

## Acknowledgements

The author wishes to thank Lindsay Childs who made several suggestions for a preliminary version of this article. The author also wishes to thank Anupam Srivastav for his assistance in originally pursuing these results. Last, the author thanks the referee for suggestions improving the presentation.